\documentclass{amsart}

\theoremstyle{definition}

\theoremstyle{remark}

\usepackage{amsfonts} \usepackage{amsmath}

\usepackage{amscd}

\begin{document}

\title[Asymptotically Efficient Estimation]{Asymptotically
 Efficient Estimation of Linear functionals in Inverse regression Models}

\author{Chris A.J. Klaassen}
\address{Korteweg-de Vries Institute for Mathematics,
University of Amsterdam, Amsterdam, THE NETHERLANDS}
\email{chrisk@science.uva.nl}
\author{Eun-Joo Lee}
\address{Department of Mathematics and Statistics,
Texas Tech University, Lubbock, TX 79409, USA}
\email{elee@math.ttu.edu}
\author{Frits H. Ruymgaart$^*$}
\address{Department of Mathematics and Statistics,
Texas Tech University, Lubbock, TX 79409, USA}
\email{ruymg@math.ttu.edu}

\date{}

\thanks{$^*$Research partly supported by National Science Foundation grant DMS-0203942.}

\begin{abstract}
\indent In this paper we will discuss a procedure to improve the usual
estimator of a linear functional of the unknown regression function in
inverse nonparametric regression models. In Klaassen \emph{et al.} (2001)
it has been proved that this traditional estimator is not asymptotically
efficient (in the sense of the H\'{a}jek - Le Cam convolution theorem)
except, possibly, when the error distribution is normal. Since this
estimator, however, is still root-n consistent a procedure in Bickel
\emph{et al.} (1993) applies to construct a modification which is
asymptotically efficient. A self-contained proof of the asymptotic
efficiency is included.
\end{abstract}

\maketitle

\setcounter{section}{0}

\section {\textbf {Introduction}}

In the nonparametric regression model with either direct or indirect
observations, both the parameter of actual interest (the unknown
regression function) and the nuisance parameter (the error density that
may also be unknown) are infinite dimensional. In this paper we consider
the problem of asymptotically efficient (AE) estimation of a linear
functional of the regression function. In other words, our aim is to
construct estimators that are asymptotically normal with smallest possible
variance. Linear functionals are of independent interest
 and often studied in the
literature (Ibragimov \& Hasminskii (1984), Goldenshluger \& Pereverzev
(2000)). They are also important because the Fourier coefficients of an
expansion of the regression function in an orthonormal basis are linear
functionals. Below we will briefly return to the latter aspect.

Suitable estimators of the regression function naturally yield estimators
of the linear functional by substitution in the inner product representing
this functional. In van Rooij \emph{et al.} (1999) it has been proved
that, when a linear functional of an indirectly sampled density is to be
estimated, substitution of a suitable density estimator produces an AE
estimator of the functional. Even when the error density is known but
arbitrary, the situation to be considered here, AE estimation of a linear
functional of the regression function turns out to be an essentially
harder problem. Klaassen \emph{et al.} (2001) have shown that substitution
of the usual type of regression function estimator does not produce an AE
estimator of the functional, except possibly when the error distribution
is normal. Another, somewhat simpler natural estimator could be easily
proposed for linear functionals, but this estimator is equivalent with and
in certain cases even identical to the plug-in estimator and hence not AE
either. It follows in particular that the usual orthonormal series type
estimator of the regression function is obtained from estimators of the
Fourier coefficients that are not AE in the above sense. Nevertheless such
estimators do in general attain the best convergence rate for the mean
integrated squared error.

In this paper we will focus on the question of how to improve the plug-in
estimator of a linear functional. The plug-in-method is a simple device to
construct an estimator of just about any functional. We will see,
moreover, that an estimator of the regression function will also be needed
in the improvement procedure. Although not AE, the plug-in estimator is
$\sqrt n$ - consistent. This means that the method to construct an AE
estimator given a $\sqrt n$ - consistent one, described in great
generality in Bickel
\emph{et al.} (1993) for infinite dimensional parameters, applies. See also Pfanzagl
(1994) for improving $\sqrt n$ - consistent estimators in models with
finite dimensional parameter. In principle our result would fit into the
general theory as described in particular in Chapter 7 of Bickel \emph{et
al.} (1993). We prefer, however, to provide a self-contained and
independent derivation in this paper. On the one hand the indirect
regression model is sufficiently specific to allow explicit calculations,
on the other hand it is of sufficient importance to warrant such an
effort. The question whether employing improved estimators of the Fourier
coefficients in a series type estimator of the regression function itself
would improve this estimator in a certain way might be of some interest.
Although this cannot be true for the rate, we suspect that improvement
will hold true at the level of constants. Investigating this matter is
beyond the scope of this paper. In Section 6, however, we will briefly
comment on this.

In order to construct an AE estimator one needs to compute the efficient
influence function. According to van der Vaart (1998) this is the
projection of the gradient of the functional onto the tangent space to the
model. See Section 3 for the details. Statement and proof that the
improved estimator is AE can be found in Section 4. Indirect nonparametric
regression occurs in a wide variety of practical situations, like
Wicksell's unfolding problem in stereology, geological prospecting,
computer tomography in imaging, just to mention a few (O'Sullivan (1986),
Kress (1989), Kirsch (1996)). Let us introduce another example (see also
Section 6.1).
\bigskip

\textbf{Example. }
Can you see the weight of a cable? To answer this question, a paraphrase
of the title of Kac's famous 1966 paper, let us first observe that the
shape of a cable suspended at its endpoints with coordinates (0,0) and
(1,0) is given by the differential equation
\begin{equation}
-\dfrac{d^2g(t)}{dt^2} = f(t),\,\,0 \leq t \leq 1,\,\,g(0) = g(1) = 0.
\tag{1.1}
\end{equation}
Apart from the sign the source term $f$ represents the load per horizontal
distance and $g$ the shape. The problem is to estimate a linear functional
$\int_0^1 f(t)\varphi(t) \: dt $, for suitable $\varphi$, like for instance
the total weight $\int_0^1 f(t) \: dt $ of the cable ($\varphi \equiv 1$);
for this special case, however, see the remark in Section 6.1.
 Estimation may be performed
by first recovering the weight distribution or source term from the data
$(X_1,Y_1),\cdots,(X_n,Y_n)$ that are independent copies of
$(X,Y)$, where
\begin{equation}
Y = g(X) + \varepsilon. \tag{1.2}
\end{equation}
We will assume that the design variable $X$ has the Uniform $(0,1)$
distribution
 and is independent of the error variable $\varepsilon$. The latter has also zero
mean, finite variance, and arbitrary density $\psi$ that is supposed to be
known.

Using the Green's function the differential equation can be rewritten in
the form
\begin{equation}
g(x) = \int_0^1 \, K(x,t)\,f(t) \: dt = (Kf)(x),\,\,0\leq x \leq 1,
\tag{1.3}
\end{equation}
where $K : L^2([0,1]) \to L^2([0,1])$ turns out to be compact and strictly
positive Hermitian with Brownian bridge covariance kernel
\begin{equation}
K(x,t) = x\wedge t-xt, \,\, 0 \leq x \leq 1,\,\, 0 \leq t \leq 1 .
\tag{1.4}
\end{equation}
The operator has strictly positive eigenvalues
\begin{equation}
\rho_k = \left \{ \dfrac{1}{\pi (k+1)} \right \}^2, \,\,k = 0,1,\cdots, \tag{1.5}
\end{equation}
with corresponding orthonormal basis of eigenfunctions
\begin{equation}
\varphi_k(x) = \sqrt 2 \sin (k+1) \pi x, \,\,0 \leq x \leq 1,\, k = 0,1,\cdots. \tag{1.6}
\end{equation}
\bigskip

The general theory below is tailored to but slightly more general than the
compact case. This generalization ensures that the noncompact identity
operator, which yields the usual direct regression, is included as a
special case.

\section{\textbf {The model, the problem, and assumptions}}

Let $(\mathbf{X},\mathcal{X},\mu)$ and $(\textbf{Z},\mathcal{Z},\nu)$ be
measurable spaces with $L^2(\mu)$ and $L^2(\nu)$ real separable
Hilbert spaces. We are given a bounded, injective linear operator $%
K:L^2(\nu)\to L^2(\mu)$. In the random design case, to be considered here,
we observe a random sample $(X_{1},Y_{1}),\cdots, (X_{n},Y_{n})$ of
independent copies of a random element $(X,Y),$ where
\begin{equation}
Y = (K f)(X)+\varepsilon = g(X)+\varepsilon, \; f \in L^2(\nu).  \tag{2.1}
\end{equation}
The indirectly observed regression or input function $f$ on
$(\textbf{Z},\mathcal{Z})$ is unknown. We will assume that
\begin{equation}
\mu (\mathbf{X})=1, \text{ and } X =_d \text{ Uniform } (\mathbf{X}). \tag{2.2}
\end{equation}
The error variable $\varepsilon$ is independent of the design $X$ with
known but arbitrary density $\psi$ with respect to Lebesgue measure. We
will assume that
\begin{equation}
\mathbf{E}\varepsilon = 0, \; \mathbf{E}\varepsilon^2 = \sigma^2 . \tag{2.3}
\end{equation}
Under these assumptions it is readily verified that $(X,Y)$ has density
\begin{equation}
p_f(x,y) = \psi(y-(Kf)(x)), \; (x,y)\in \mathbf X \times \mathbb R , \; f
\in  L^2(\nu). \tag{2.4}
\end{equation}

A star attached to an operator denotes its adjoint. The operator
$R=(K^*K)^{\frac{1}{2}} : L^2(\nu)\rightarrow L^2(\nu)$, then, is strictly
positive Hermitian. We will assume that there exists an orthonormal basis
for $L^2(\nu)$ consisting of eigenfunctions
$\varphi_0,\varphi_1,\cdots$ of the operator $R$ with corresponding eigenvalues
$\rho_0, \rho_1,\cdots > 0$ satisfying
\begin{equation}
\text{sup}_{k \ge 0} \rho_k < \infty . \tag{2.5}
\end{equation}
On the one hand compact operators $R$ are included, since for those a
basis exists with eigenvalues satisfying $\rho_k \downarrow 0 , \text{ as
} k \to \infty $ (Debnath \& Mikusi\'nski (1999)). On the other hand the
direct model with $K = I$ and hence $R = I$, where
$I$ is the identity operator, is included as well. The operator
$I$ is not compact but satisfies the condition above for any basis, with
$\rho_k = 1$ for all $k$. According to the polar decomposition
(Riesz \& Nagy (1990)) there exists a partial isometry
$V : L^2(\nu) \to L^2(\mu)$ such that $K = VR$ and $K^* = RV^*$.
It should be noted that $V^*V$ is the identity on the range of $R$.
Let us write
\begin{equation}
\varphi_{V,k} = V\varphi_k , \tag{2.6}
\end{equation}
and observe that the $\varphi_{V,k}$ are orthonormal in $L^2(\mu)$.

The problem to be considered here is estimation of a linear functional
$ f \mapsto \langle f , \varphi \rangle , f \in L^2(\nu) ,$ for some given
$\varphi \in L^2(\nu)$ with $\|\varphi\| = 1.$
In view of $\langle f , \varphi
\rangle $ $= \sum_{k \ge 0} \langle f , \varphi_k \rangle
\langle \varphi , \varphi_k \rangle$ it seems plausible that estimation of the
special linear functionals defined by the Fourier coefficients
\begin{equation}
f_k = \langle f , \varphi_k \rangle  \tag{2.7}
\end{equation}
might suffice. This is in fact true under an
extra condition, and some details can be found in Section 5. Hence we will
focus on estimating an arbitrary Fourier coefficient. As a generic example
let us consider the functional
\begin{equation}
f \mapsto f_0 = \langle f , \varphi _0 \rangle, f \in L^2(\nu),
\tag{2.8}
\end{equation}
where $\varphi_0$ is the first basis element. It is useful to observe
\begin{equation}
f_k = \langle f , \varphi _k \rangle = \langle f, K^*VR^{-1}\varphi _k
\rangle = \langle Kf, VR^{-1} \varphi _k \rangle
= \dfrac{1}{\rho _k} \langle Kf, \varphi _{V,k} \rangle . \tag{2.9}
\end{equation}

The assumptions below will be briefly discussed in Section 6. The basis
elements are supposed to satisfy the uniform boundedness conditions
\begin{equation}
\text{sup}_{z \in \mathbf Z, k \geq 0}
 |\varphi_k(z)| < \infty \tag{2.10}
\end{equation}
and
\begin{equation}
\text{sup}_{x \in \mathbf X, k \geq 0}
 |\varphi_{V,k}(x)| < \infty. \tag{2.11}
\end{equation}

Regarding the input function $f$ it will be assumed that there exists a
sequence $(m(n))_{n\ge 1}$ satisfying
\begin{equation}
m = m(n) \to \infty \; \text{ and }\; \dfrac{m}{\sqrt n} \to 0 , \text{ as
} n
\to
\infty ,
\tag{2.12}
\end{equation}
for which
\begin{equation}
\sqrt n \sum_{k>m} f_k ^2 \to 0 , \text{ as } n \to \infty , \tag{2.13}
\end{equation}
\begin{equation}
\sum_{k>m} |f_k| \to 0 , \text{ as } n \to \infty . \tag{2.14}
\end{equation}

The error density $\psi$ is supposed to be twice differentiable, i.e.,
\begin{equation}
\psi^{''} \text{ exists on } \mathbb R .
\tag{2.15}
\end{equation}
Denoting the score function for location by
\begin{equation}
\Lambda = -\dfrac{\psi '}{\psi}, \tag{2.16}
\end{equation}
we will also need that
\begin{equation}
\Lambda ' \text{ and } \Lambda^{''} \text{ exist and are bounded on } \mathbb R , \tag{2.17}
\end{equation}
with finite Fisher information
\begin{equation}
\mathbf E \Lambda ^2 (\varepsilon) = \mathbf E \Lambda ' (\varepsilon) =
\int_{-\infty}^{\infty} \left (\dfrac{\psi '}{\psi} \right )^2(y)
\psi(y) \: dy = I(\psi) < \infty . \tag{2.18}
\end{equation}
For the first equality in (2.18) we need (2.15).

\section{\textbf {Construction of an asymptotically efficient estimator}}

The input function $f \in L^2(\nu)$ has the $L^2$ - expansion
\begin{equation}
f(z) = \sum_{k \ge 0} f_k \varphi_k(z), \; z \in \mathbf Z , \tag{3.1}
\end{equation}
in the orthonormal basis of eigenfunctions, and the usual estimator of the
regression function in this context is given by
\begin{equation}
\hat{f}(z) = \hat{f}_{(m)}(z) = \sum_{k \le m} \hat{f_k}\varphi_k(z), \; z \in
\mathbf Z , \tag{3.2}
\end{equation}
for suitable $m = m(n) \to \infty$, as $n \to \infty$, where $\hat{f_k}$
is the estimator of $f_k$ (see (2.1) and (2.9)) given by
\begin{equation}
\hat{f_k} = \dfrac{1}{n} \sum_{i=1}^n \dfrac{1}{\rho_k} Y_i
\varphi_{V,k}(X_i).  \tag{3.3}
\end{equation}
See, for instance, Johnstone \& Silverman (1990). Since we are estimating
$f_0$ in (2.8), a plug-in estimator simply equals
\begin{equation}
\langle \hat{f} , \varphi_0 \rangle = \hat{f_0} , \tag{3.4}
\end{equation}
as given by (3.3) for $ k = 0 $.

The estimators $\hat{f_k} ( k = 0,1,\cdots ) $ have some desirable
properties. Because of
\begin{align}
\mathbf E \hat{f_k} &= \dfrac{1}{\rho_k} \mathbf E Y\varphi_{V,k}(X)
                     = \dfrac{1}{\rho_k} \mathbf E (Kf)(X)\varphi_{V,k}(X) \tag{3.5}\\
                    &= \dfrac{1}{\rho_k} \langle Kf , V\varphi_k\rangle
                     = \dfrac{1}{\rho_k} \langle f , RV^*V\varphi_k \rangle \notag \\
                    &= \dfrac{1}{\rho_k} \rho_k \langle f, \varphi_k \rangle
                     = f_k  \notag
\end{align}
they are unbiased, but $\hat{f} = \hat{f}_{(m)}$ is not; let us write
\begin{equation}
f_{(m)}(z) = \mathbf E \hat{f}_{(m)}(z) = \sum_{k \le m} f_k
\varphi_k(z), \; z \in \mathbf Z . \tag{3.6}
\end{equation}
Furthermore we have
\begin{align}
&\mathbf E (\hat{f_k} - f_k)^2 = \textbf{ Var }\hat{f_k} =
\dfrac{1}{n} \textbf{ Var } \dfrac{1}{\rho_k} Y \varphi_{V,k}(X)
\tag{3.7} \\
&\le \dfrac{1}{n} \left ( \dfrac{1}{\rho_k} \right )^2 \mathbf E Y^2
\varphi_{V,k}^2 (X)  \le C \dfrac{1}{n} \left ( \dfrac{1}{\rho_k} \right
)^2 , \notag \end{align} where $0 < C < \infty$ will throughout be used as
a generic constant that will not depend on $n$ or $k$ and that here can be
taken equal to
\begin{equation}
C = (\text{ sup }_{x \in \mathbf X, k \ge 0} |\varphi_{V,k}^2 (x)|)
\cdot (\|Kf\|^2 + \sigma^2), \tag{3.8}
\end{equation}
by assumption (2.11).

The central limit theorem yields at once the asymptotic normality of the
empirical Fourier coefficients. In particular we have
\begin{equation}
\sqrt n ( \hat{f_0} - f_0 ) \to _d \text { Normal } (0,
\sigma_0^2(f)), \text{ as } n \to \infty , \tag{3.9}
\end{equation}
where (cf. (2.1), (2.3), and (2.6))
\begin{align}
\sigma_0^2(f) &= \left ( \dfrac{1}{\rho_0} \right )^2
\textbf { Var} \left( Y \varphi_{V,0}(X)\right) \tag{3.10}\\
&\geq \left ( \dfrac{1}{\rho_0} \right )^2
\mathbf E \textbf { Var} \left( Y \varphi_{V,0}(X)|X\right) \notag \\
&= \left ( \dfrac{\sigma}{\rho_0} \right )^2 \
\mathbf E \,\varphi_{V,0}^2 (X)
= \left ( \dfrac{\sigma}{\rho_0} \right )^2  \notag
\end{align}
holds with a strict inequality, unless
$\textbf { Var} \left(\mathbf E \left( Y
\varphi_{V,0}(X)|X\right)\right)$ vanishes, i.e., unless
$(Kf)(X)\varphi_{V,0}(X)$ is degenerate.

It has already been observed in Klaassen
\emph{et al.} (2001) that
$\sigma_0^2(f)$ is in general strictly larger than the optimal variance
according to the H\'{a}jek - Le Cam convolution theorem (van der Vaart
(1998)). In other words the estimator $\hat{f_0}$ is not asymptotically
efficient, but it is $\sqrt n$ - consistent. Such estimators can be
improved. In order to do so we first need to briefly review some results
from Klaassen
\emph{et al.} (2001).

Since we assume $\psi$ to be known, in terms of square roots of the
densities the model is
$\mathcal{S} = \{s_f, f \in L^2(\nu) \}$, with $s_f=\sqrt{p_f}$ and
$p_f$ as in (2.4). The tangent space at
$f \in L^2(\nu)$ to this model is given by
\begin{equation}
\overset{\bullet }{\mathcal S}_f = [ \overset{\bullet }s_{f,h} ,\;
h \in \mathfrak{R}_K ] = \{ \overset{\bullet }s_{f,h} ,\; h \in
\overline{\mathfrak{R}_K}\}
\subset L^2(\mu \times \lambda) ,    \tag{3.11}
\end{equation}
where $\lambda$ is Lebesgue measure on $\mathbb{R}$, $\mathfrak{R}_K$ is
the range of $K$,
$\overline{\mathfrak{R}_K}$ its closure, and
\begin{equation}
\overset{\bullet }s_{f,h}(x,y) = -\dfrac{\psi ' (y-(Kf)(x))}{2s_f(x,y)}h(x),
\; x \in \mathbf X,\, y \in \mathbb R . \tag{3.12}
\end{equation}
The gradient of the functional $f \mapsto T(s_f) = \langle f ,
\varphi_0 \rangle$ at
$f$ is given by (cf. (2.9))
\begin{equation}
\overset{\bullet}{T}_f(x,y) = \dfrac{2}{\rho_0}(y-(Kf)(x))
\varphi_{V,0}(x)s_f(x,y), \; x \in \mathbf X , \, y \in  \mathbb R . \tag{3.13}
\end{equation}
Let $\widetilde{T}_f$ denote the projection in $L^2(\mu \times \lambda) $
of
 $\, \overset{\bullet}{T}_f$ onto $\overset{\bullet }{\mathcal S}_f$. Then
the optimal variance mentioned in the preceding paragraph equals
\begin{equation}
\widetilde{\sigma}_0^2 (f) = \dfrac{1}{4} \| \widetilde{T}_f \|^2 .
\tag{3.14}
\end{equation}
For these results see Klaassen \emph{et al.} (2001).

In order to construct an estimator with limiting normal distribution
having the variance in (3.14) we first need to explicitly compute
$\widetilde{T}_f$. Because
$\widetilde{T}_f$ has to belong to $\overset{\bullet }{\mathcal S}_f$ it is of the form
\begin{equation}
\widetilde{T}_f(x,y) = -\dfrac{\psi ' (y-(Kf)(x))}{2s_f(x,y)}\widetilde h(x),
\; x \in \mathbf X,\, y \in \mathbb R , \tag{3.15}
\end{equation}
for some $\widetilde h \in \overline{\mathfrak{R}_K}$, and because of
$\overset{\bullet}{T}_f -
\widetilde{T}_f \perp \overset{\bullet }{\mathcal S}_f$ we must have (see
also (3.12))
\begin{align}
\langle \overset{\bullet}{T}_f - \widetilde{T}_f , \overset{\bullet }s_{f,h} \rangle
&=\iint \{ \dfrac{2}{\rho_0} (y-(Kf)(x))\varphi_{V,0}(x) s_f(x,y) +
\dfrac{\psi '
(y-(Kf)(x))}{2s_f(x,y)} \widetilde h(x) \} \tag{3.16} \\
&\qquad \times \{ \dfrac{\psi ' (y-(Kf)(x))}{2s_f(x,y)} h(x) \}
\,\, d\mu(x)\, dy = 0 , \notag
\end{align}
for all  $h \in \overline{\mathfrak{R}_K}$. It follows from (2.3) and
(2.15) that
$\int_{-\infty}^\infty y\psi ' (y) \: dy = -1$ holds; cf. Lemma I.2.4.b of
H\'ajek and \u{S}id\'ak (1967). Exploiting this fact and using the
notation
$I(\psi)$ defined in (2.18) straightforward integration shows that (3.16)
entails
\begin{equation}
\langle I(\psi)\widetilde h - \dfrac{4}{\rho_0}\varphi_{V,0} , h \rangle = 0 ,
\text{ for all } h \in
\overline{\mathfrak{R}_K} . \tag{3.17}
\end{equation}
It is obvious that $\dfrac{1}{\rho_0}\varphi_{V,0} \in \mathfrak{R}_K$.
Since $\widetilde h \in \overline{\mathfrak{R}_K}$ it follows from (3.17)
that
\begin{equation}
\widetilde h(x) = \dfrac{4}{I(\psi)\rho_0}\varphi_{V,0}(x) ,
\; x \in \mathbf X . \tag{3.18}
\end{equation}
Combination with (3.14) yields
\begin{equation}
\widetilde{\sigma}_0^2 (f) = \dfrac{1}{4} \| \widetilde{T}_f \|^2 =
\dfrac{1}{\rho_0^2 I(\psi)} , \tag{3.19}
\end{equation}
for the value of the optimal variance. Note that
\begin{equation}
1=\left(\int_{-\infty}^\infty y\psi ' (y) \: dy \right)^2 \leq
\sigma^2I(\psi) \notag
\end{equation}
holds by the Cram\'er - Rao inequality with equality if and only if
$\psi ' (Y)/\psi(Y)$ is linear in $Y$ a.s. under $\psi$, i.e. if and only
if $\psi$ is a normal density.

Consequently, it is immediate from (3.10) that the string of inequalities
\begin{equation}
\sigma_0^2(f) \ge \left (\dfrac{1}{\rho_0} \right )^2 \sigma^2
\ge \left (\dfrac{1}{\rho_0} \right )^2 \dfrac{1}{I(\psi)} = \widetilde
\sigma_0^2(f)  \tag{3.20}
\end{equation}
holds. For nonnormal densities $\psi$ the second inequality is strict, and
the estimator
$\hat{f_0}$ then turns out not to be asymptotically efficient. If $\psi$
is normal, we have $I(\psi) = 1 / \sigma^2$, so that the second inequality
in (3.20) is an equality. However, asymptotic efficiency of
$\hat{f_0}$ remains impossible, since the first inequality in
(3.20) cannot be an equality for all $f$, as argued in (3.10). See
Klaassen \emph{et al.} (2001) for some details.

We are now in a position to construct the improved estimator. According to
van der Vaart (1998) define the efficient influence function by
\begin{align}
\widetilde{T}_f(x,y) \cdot \dfrac{1}{2 s_f(x,y)} &=
-\dfrac{1}{I(\psi)} \dfrac{\psi ' (y-(Kf)(x))}{\psi
(y-(Kf)(x))}\dfrac{1}{\rho_0}\varphi_{V,0}(x)
\tag{3.21}  \\
&= \dfrac{1}{I(\psi)} \Lambda (y-(Kf)(x))\dfrac{1}{\rho_0}\varphi_{V,0}(x)
, \; x \in \mathbf X , y \in \mathbb R , \notag
\end{align}
using the notation $\Lambda$ introduced in (2.16). Following a procedure
in Bickel \emph{et al.} (1993, Chapter 7) let us now introduce the
estimator
\begin{equation}
\Hat{\Hat{f_0}} = \hat{f_0} +
\dfrac{1}{n} \sum_{i=1}^n  \dfrac{1}{\rho_0 I(\psi)} \Lambda (Y_i -(K \hat
f_{(m)})(X_i)) \varphi_{V,0}(X_i), \tag{3.22}
\end{equation}
where $\hat{f}_{(m)} = \hat{f}$ is defined in (3.2). We will see in the
next section that this turns out to be an asymptotically efficient
estimator of
$f_0$.

\section{\textbf {The main theorem}}

\textbf{Theorem. } \emph{Suppose that all the assumptions listed in Section 2 are
fulfilled. Then }$\Hat{\Hat{f_0}}$ \emph{defined in} (3.22) \emph{is an
asymptotically efficient estimator of} $f_0$, \emph{i.e.}
\begin{equation}
\sqrt n ( \Hat{\Hat{f_0}} - f_0 ) \to_d \text{ Normal } \left( 0 ,
\dfrac{1}{\rho_0^2 I(\psi)} \right) , \text{ as } n \to \infty , \tag{4.1}
\end{equation} \emph{where the variance in the normal distribution is optimal.}

\textbf{Proof. }
A Taylor expansion yields
\begin{equation}
\sqrt n ( \Hat{\Hat{f_0}} - f_0 ) = A_n + Q_n + R_n  \tag{4.2}
\end{equation}
with
\begin{equation}
A_n = \dfrac{1}{\sqrt n} \sum_{i=1}^n \dfrac{1}{\rho_0 I(\psi)} \Lambda
(\varepsilon_i)  \varphi_{V,0} (X_i), \tag{4.3}
\end{equation}
\begin{align}
Q_n &= \sqrt n [ \langle \hat{f}-f , \varphi_0 \rangle \tag{4.4}\\
& \quad - \dfrac{1}{n} \sum_{i=1}^n
 \dfrac{1}{\rho_0 I(\psi)} \{ (K\hat f)(X_i) - (Kf)(X_i) \} \Lambda '
(\varepsilon_i) \varphi_{V,0} (X_i) ] , \notag
\end{align}
\begin{equation}
R_n = \dfrac{1}{2\sqrt n} \sum_{i=1}^n
 \dfrac{1}{\rho_0 I(\psi)} \{ (K\hat f)(X_i) - (Kf)(X_i) \}^2 \Lambda ''
(\widetilde{\varepsilon_i})  \varphi_{V,0} (X_i) . \tag{4.5}
\end{equation}
In (4.5), $\widetilde{\varepsilon_i}$ is a random variable between $Y_i -
(Kf)(X_i) = \varepsilon_i$ and $Y_i - (K\hat f)(X_i)$.

It follows from (2.11) that
\begin{align}
\text{sup}_{x \in \mathbf X} | \varphi_{V,0}(x) | < \infty . \tag{4.6}
\end{align}
Together with (2.17) we find ( $C$ generic ! )
\begin{align}
|R_n| &\le \dfrac{C}{\sqrt n} \sum_{i=1}^n \{(K \hat f)(X_i) - (Kf)(X_i)\}^2 \tag{4.7} \\
&= \dfrac{C}{\sqrt n} \sum_{i=1}^n \{ \sum_{k \le m} (\hat{f_k} -
f_k)\rho_k \varphi_{V,k}(X_i)
- \sum_{k > m} f_k \rho_k \varphi_{V,k}(X_i) \}^2  \notag \\
&\le \dfrac{C}{\sqrt n} \sum_{k \le m} \sum_{l \le m} (\hat{f_k} - f_k) (\hat{f_l} - f_l)
\rho_k \rho_l \{\sum_{i=1}^n \varphi_{V,k}(X_i) \varphi_{V,l}(X_i)\} \notag \\
& \qquad +\dfrac{C}{\sqrt n} \sum_{k > m} \sum_{l > m} f_k f_l \rho_k \rho_l
\{\sum_{i=1}^n \varphi_{V,k}(X_i) \varphi_{V,l}(X_i)\} \notag \\
&= R_{n1} + R_{n2} . \notag
\end{align}
It clearly suffices to show that $\mathbf E R_{nj} \to 0$, as $n \to
\infty$, for
$j = 1,2,$ in order to obtain $R_n \to_p 0$, as $n \to \infty$.

Regarding the first term we see that
\begin{align}
\mathbf E R_{n1} &= \dfrac{C}{n^{5/2}} \sum_{k \le m} \sum_{l \le m} [
\sum_{\alpha = 1}^n \sum_{\beta = 1}^n \sum_{i=1}^n
\mathbf E \{Y_\alpha \varphi_{V,k}(X_\alpha) - f_k \rho_k \} \tag{4.8} \\
& \qquad \times \{Y_\beta \varphi_{V,l}(X_\beta) - f_l \rho_l \}
\{\varphi_{V,k}(X_i) \varphi_{V,l}(X_i) \}] . \notag
\end{align}
By decomposing $\sum \sum \sum_{\alpha, \beta, i} = \sum \sum \sum_{\alpha
\ne \beta
\ne i} + \sum \sum_{\alpha \ne \beta = i} + \sum \sum_{\beta \ne \alpha = i}+\\
\sum \sum_{\alpha = \beta \ne i} + \sum_{\alpha = \beta = i}$, and by realizing
that the random variables labeled with $\alpha$ and $\beta$ are centered at
$0$ and that the $\varphi_{V,k}$ are orthonormal, we arrive at ($\delta_{kl}$ is
Kronecker's delta)
\begin{align}
0 &\le \mathbf E R_{n1} \tag{4.9}\\
&= \dfrac{C}{n^2 \sqrt n} \:\sum_{k \le m} \sum_{l \le m}
\{ n(n-1) \delta_{kl} \textbf{Cov}(Y\varphi_{V,k}(X),Y\varphi_{V,l}(X)) \notag\\
& \qquad + n \mathbf E (Y\varphi_{V,k}(X) - f_k \rho_k) (Y\varphi_{V,l}(X)
- f_l \rho_l) \varphi_{V,k}(X) \varphi_{V,l}(X)\} \notag \\
&\le \dfrac{C}{n^2 \sqrt n} [ n(n-1) \sum_{k \le m} \mathbf E Y^2 \varphi_{V,k}^2 (X) \notag \\
& \qquad + n \sum_{k \le m} \sum_{l \le m} \{ \mathbf E Y^2 \varphi_{V,k}^2 (X)
\mathbf E Y^2 \varphi_{V,l}^2 (X) \}^{1/2} ] \notag \\
&\le C \left( \dfrac{m}{\sqrt n} + \dfrac{m^2}{n\sqrt n} \right ) \to 0 ,
\text{ as } n \to \infty . \notag
\end{align}
Here we have used (2.11) and (2.12). See also the calculation in (3.7).


Next let us observe that
\begin{align}
0 \le \mathbf E R_{n2} &= \dfrac{C}{\sqrt n} \: \sum_{k > m} \sum_{l > m}
n f_k f_l \rho_k \rho_l
\delta_{kl} \tag{4.10} \\
&\le C \sqrt n \sum_{k > m} f_k^2 \rho_k^2  \notag \\
&\le C \sqrt n \sum_{k > m} f_k^2 \to 0 , \text{ as } n \to \infty ,
\notag
\end{align}
by assumption (2.13). This settles the asymptotic negligibility of $R_n$.

For brevity let us introduce
\begin{equation}
\overline{U}_k = \dfrac{1}{n} \sum_{i=1}^n \{ \langle \varphi_0 , \varphi_k \rangle -
\dfrac{1}{\rho_0 I(\psi)} \rho_k \varphi_{V,k}(X_i) \Lambda ' (\varepsilon_i)
\varphi_{V,0} (X_i) \}, \tag{4.11}
\end{equation}
and note that $\overline{U}_k$ is an average of i.i.d. random variables
with zero mean. Indeed we have (cf.(2.18))
\begin{align}
&\langle \varphi_0 , \varphi_k \rangle -
\dfrac{1}{\rho_0 I(\psi)} \rho_k \mathbf E \varphi_{V,k}(X) \Lambda '
(\varepsilon) \varphi_{V,0}(X) \tag{4.12}\\
&= \delta_{0k} - \dfrac{1}{\rho_0 I(\psi)} \rho_k I(\psi) \delta_{k0} = 0
,
\text{ for all }k.
\notag
\end{align}
Because of
\begin{align}
Q_n &= \sqrt n \sum_{k \le m} (\hat{f_k} - f_k) \overline{U}_k -
\sqrt{n} \sum_{k > m} f_k \overline{U}_k  \tag{4.13} \\
&= Q_{n1} + Q_{n2} , \notag
\end{align}
in order to show that $Q_n \to_p 0$, as $n \to \infty$, it suffices to
prove $\mathbf E |Q_{n1}| \to 0$ and $\mathbf E |Q_{n2}| \to 0$.

In order to deal with $Q_{n1}$ let us first note that
\begin{align}
&\mathbf E (\hat{f_k} - f_k)^2 \, \overline{U}_k^2 = \tag{4.14}\\
&= \dfrac{1}{n^4} \sum_{i=1}^n \sum_{j=1}^n \sum_{\alpha=1}^n
\sum_{\beta=1}^n
\mathbf E \{ Y_i \varphi_{V,k}(X_i) - \rho_k f_k \}
\{ Y_j \varphi_{V,k}(X_j) - \rho_k f_k \} \notag \\
& \qquad \times
 \{ \langle \varphi_0 , \varphi_k \rangle /{\rho_k} - \dfrac{1}{\rho_0 I(\psi)}
\varphi_{V,k} (X_\alpha) \Lambda ' (\varepsilon_\alpha) \varphi_{V,0}
(X_\alpha) \} \notag \\
& \qquad \times \{ \langle \varphi_0 , \varphi_k \rangle / {\rho_k} -
\dfrac{1}{\rho_0 I(\psi)} \varphi_{V,k} (X_\beta) \Lambda ' (\varepsilon_\beta)
\varphi_{V,0}(X_\beta) \} . \notag
\end{align}
Since each factor has zero expectation the only terms that contribute are
those with all 4 indices equal and those with 2 pairs of indices (but not
all 4) equal. By (2.11), (2.17), and (3.7) the contribution of the first
group is seen to be bounded by $C \cdot n$ and the second group is seen to
yield a contribution bounded by $C \cdot n^2$. This entails
\begin{equation}
\mathbf E (\hat{f_k} - f_k)^2 \,
\overline{U}_k^2  \le \dfrac{1}{n^4} (C \cdot n + C \cdot n^2) \le
\dfrac{C}{n^2} , \tag{4.15}
\end{equation}
and hence, by applying Schwarz's inequality,
\begin{align} \mathbf E |Q_{n1}| &\le {\sqrt n} \sum_{k \le m} \{
\mathbf E (\hat{f_k} - f_k)^2 \, \overline{U}_k^2 \}^{1/2} \tag{4.16}\\
&\le {\sqrt n} m \dfrac{C}{n} = C \dfrac{m}{\sqrt n} \to 0,
\text{ as } n \to \infty , \notag
\end{align}
by assumption (2.12).

For $Q_{n2}$ observe that,
\begin{align}
\mathbf E \overline{U}_k^2 &= \dfrac{1}{n^2} \sum_{i=1}^n \sum_{j=1}^n
\mathbf E \{ \langle \varphi_0 , \varphi_k \rangle -
\dfrac{1}{\rho_0 I(\psi)} \rho_k \varphi_{V,k}(X_i) \Lambda ' (\varepsilon_i)
\varphi_{V,0} (X_i) \} \tag{4.17} \\
& \qquad \times\{ \langle \varphi_0 , \varphi_k \rangle -
\dfrac{1}{\rho_0 I(\psi)} \rho_k \varphi_{V,k}(X_j) \Lambda ' (\varepsilon_j)
\varphi_{V,0}(X_j) \} \notag \\
&\le \dfrac{C}{n} \left (\dfrac{\rho_k}{\rho_0} \right )^2 \mathbf E
\varphi_{V,k}^2(X) (\Lambda ')^2 (\varepsilon)
\varphi_{V,0}^2(X) \notag \\
&\le \dfrac{C}{n} \rho_k^2 . \notag
\end{align}
This yields, again by the Schwarz inequality
\begin{align}
\mathbf E |Q_{n2}| &\le {\sqrt n} \sum_{k > m} |f_k|
\{ \mathbf E \overline{U}_k^2 \}^{1/2} \tag{4.18}\\
&\le C \sum_{k > m} |f_k| \rho_k \to 0, \text{ as } n \to \infty , \notag
\end{align}
by assumptions (2.5) and (2.14).

Finally let us consider $A_n$ in (4.3). Since the terms are i.i.d. with
\begin{equation}
\mathbf E \dfrac{1}{\rho_0 I(\psi)} \Lambda (\varepsilon) \varphi_{V,0}(X) =
0 , \tag{4.19}
\end{equation}
\begin{align}
\textbf {Var} \dfrac{1}{\rho_0 I(\psi)} \Lambda (\varepsilon) \varphi_{V,0}(X)
&= \left (\dfrac{1}{\rho_0 I(\psi)} \right )^2 \mathbf E \Lambda^2
(\varepsilon) \mathbf E \varphi_{V,0}^2(X) \tag{4.20} \\
&= \dfrac{1}{\rho_0^2 I(\psi)} , \notag
\end{align}
the central limit theorem entails at once that
\begin{equation}
A_n \to_d \text{ Normal } \left( 0 , \dfrac{1}{\rho_0^2 I(\psi)} \right) ,
\text{ as } n \to \infty . \tag{4.21}
\end{equation}
Because we have seen that $R_n \to_p 0$ and $Q_n \to_p 0$, this is also
the limiting distribution of the expression on the left in (4.2), as was
to be shown.

\section{\textbf {Estimating an arbitrary linear functional}}

In this section we want to consider the problem of estimating the linear
functional $f_\varphi = \langle f , \varphi \rangle$, for an arbitrary
$\varphi \in L^2(\nu)$ with $\|\varphi\| = 1$. Since $f_\varphi = \sum_{k \ge
0} \langle f , \varphi_k \rangle  \langle \varphi , \varphi_k \rangle =
\sum_{k \ge 0} f_k \langle \varphi , \varphi_k \rangle$,
in view of the preceding results we expect
\begin{equation}
\Hat{\Hat{f}}_\varphi = \sum_{k \ge 0} \Hat{\Hat{f}}_k \langle \varphi ,
\varphi_k \rangle , \tag{5.1}
\end{equation}
to be an asymptotically efficient estimator. Before proceeding we need to
introduce the extra condition that
\begin{equation}
\sum_{k \ge 0} \dfrac{1}{\rho_k} | \langle \varphi , \varphi_k \rangle | <
\infty . \tag{5.2}
\end{equation}
This condition ensures that $\varphi$ is in the range of $(K^{-1})^* .$
Let us write
\begin{align}
\gamma &= (K^{-1})^* \varphi = VR^{-1} (\sum_{k \ge 0}
\langle \varphi , \varphi_k \rangle \varphi_k) \tag{5.3}\\
&= \sum_{k \ge 0} \dfrac{1}{\rho_k}  \langle \varphi , \varphi_k \rangle
\varphi_{V,k}  .
\notag
\end{align}
Thanks to (2.10) and (5.2) the convergence in (5.3) is even pointwise.

To verify the asymptotic efficiency let us first observe that the optimal
variance in the normal component of regular sequences of estimators equals
\begin{equation}
\widetilde{\sigma}_\varphi^2 (f) = \dfrac{\| \gamma \|^2}{I(\psi)} . \tag{5.4}
\end{equation}
This can be shown by virtually the same method as employed in Section 3
for
$\varphi_0$.

Writing the decomposition in (4.2) for general $\varphi_k$ (rather than
for
$\varphi_0$) as
\begin{equation}
\sqrt n ( \Hat{\Hat{f}}_k - f_k ) = A_{n,k} + Q_{n,k} + R_{n,k} ,
\tag{5.5}
\end{equation}
we have
\begin{align}
\sqrt n ( \Hat{\Hat{f}}_\varphi - f_\varphi )&= \sum_{k \ge 0}
\langle \varphi , \varphi_k \rangle A_{n,k} + \sum_{k \ge 0}
\langle \varphi , \varphi_k \rangle (Q_{n,k} + R_{n,k}) \tag{5.6} \\
&= S_{n,1} + S_{n,2} . \notag
\end{align}
As follows from (5.3),
\begin{align}
S_{n,1} &= \dfrac{1}{\sqrt n} \sum_{i=1}^n \dfrac{1}{I(\psi)} \Lambda
(\varepsilon_i) \sum_{k \ge 0} \dfrac{1}{\rho_k}
\langle \varphi , \varphi_k \rangle \varphi_{V,k}(X_i) \tag{5.7} \\
&= \dfrac{1}{\sqrt n} \sum_{i=1}^n \dfrac{1}{I(\psi)} \Lambda
(\varepsilon_i) \gamma(X_i) , \notag
\end{align}
which means that $S_{n,1}$ is well defined. Exploiting the calculations in
the proof of the Theorem in Section 4 it follows that $S_{n,2}$ is also
well defined and that
\begin{equation}
\sum_{k \ge 0} \langle \varphi , \varphi_k \rangle ( Q_{n,k} + R_{n,k} ) =
o_p(1) . \tag{5.8}
\end{equation}
Combination yields the following results.

\textbf{Theorem. } \emph{Suppose that condition (5.2) is fulfilled in addition
to the assumptions in Section 2.  Then }$\Hat{\Hat{f}}_\varphi$ \emph{is
an asymptotically efficient estimator of} $f_\varphi$, \emph{i.e.}
\begin{equation}
\sqrt n ( \Hat{\Hat{f}}_\varphi - f_\varphi ) \to_d \text{
Normal } \left( 0 , \dfrac{\|  \gamma \|^2}{I(\psi)} \right) , \text{ as }
n \to \infty , \tag{5.9}
\end{equation}

\section{\textbf {Some comments on the conditions and \\
improved regression estimation}}

\textbf{6.1. The conditions.} In the example of Section 1 the
operator $K$ itself is Hermitian and hence $V = I$ so that
\begin{equation}
\varphi_k(x) = \varphi_{V,k}(x) = \sqrt 2 \sin (k+1) \pi x ,\, 0 \le x \le 1 ,\, k
= 0,1,\cdots \tag{6.1}
\end{equation}
It is immediate that conditions (2.5), (2.10), and (2.11) are fulfilled.
It should be noted, however, that the function $\varphi \equiv 1$ on
$[0,1]$ doesn't satisfy (5.2). If the total weight is to be estimated, one
should therefore employ a sufficiently smooth (near 0 and 1) approximation
of this function.

Next suppose that the input function $f$ satisfies
$|\langle f , \varphi_k \rangle | = |f_k| \asymp k^{-s}$ for some
$ s > \frac{1}{2} ,$ and that $m(n) \asymp n^r$ for some $0 < r < \frac{1}{2}
$. Then we have
\begin{equation}
\sqrt n \sum_{k>m} f_k^2 = O(n^{1/2 + r(1-2s)}), \tag{6.2}
\end{equation}
\begin{equation}
\sum_{k>m} |f_k| = O(n^{r(1-s)}), \; s > 1 , \tag{6.3}
\end{equation}
and apparently conditions (2.13) and (2.14) are satisfied if
\begin{equation}
0 < r < \dfrac{1}{2} \,, \;\; s > \dfrac{1+2r}{4r} \,. \tag{6.4}
\end{equation}

The conditions on the error density also appear to be usually fulfilled. A
nonnormal density that satisfies conditions (2.15) - (2.18) is, for
instance, the logistic density
\begin{equation}
\psi (x) = \dfrac{e^{-x} }{ (1 + e^{-x})^2 }, \; x \in \mathbb R . \tag{6.5}
\end{equation}
In particular $\Lambda ' \text{ and } \Lambda^{''}$ turn out to be bounded
indeed.

\textbf{6.2. Improving regression estimation.}
It should be noted that for the present results it is irrelevant whether
the input estimator $\hat f_{(m)}$ attains the optimal MISE rate. It is
not hard to see, however, that conditions on $m$ and $f$ ensuring this
rate to be optimal are in general compatible with those in (2.12) -
(2.14). In this discussion, however, we will allow the truncation index $M
= M(n)$ of the traditional input estimator $\hat f_{(M)}$, that attains
the optimal MISE rate, to tend to infinity at a different rate than the $m
= m(n)$ used above.

The results in this paper regard the variances of the limiting normal
distributions of the estimators, and not the variances or MSE's of these
estimators themselves. It is clear, however, that
\begin{equation}
\text{MSE}(\hat f_{0}) = \textbf{Var}\hat f_{0} = \dfrac{1}{n}
\sigma_0^2(f), \tag{6.6}
\end{equation}
and we conjecture that
\begin{equation}
\text{MSE}(\Hat{\Hat{f_0}}) = \mathbf{E} (\Hat{\Hat{f_0}} - f_0)^2 =
\dfrac{1}{n} \dfrac{1}{I(\psi)} + o(\dfrac{1}{n}) . \tag{6.7}
\end{equation}
This would imply that $\Hat{\Hat{f_0}}$ improves $\hat f_0$ also with
respect to the MSE (at the level of constants). We similarly expect each
of the $\Hat{\Hat{f}}_k$ to improve $\hat f_k$ regarding MSE, and
eventually
\begin{equation}
\Hat{\Hat{f}}_{(M)}(y) =\sum_{k \le M} \Hat{\Hat{f}}_k \varphi_k (y), y \in
\mathbf{Y}, \tag{6.8}
\end{equation}
to improve $\hat f_{(M)}$ with respect to MISE (at the level of
constants). The actual calculations leading to (6.7) will differ from
those in Section 4 and might be lengthy. Moreover, if theoretically
$\Hat{\Hat{f}}_{(M)}$ would improve $\hat f_{(M)}$, it would be
interesting to perform simulations for several nonnormal error
distributions, to get an insight into the difference of the performance
for finite sample sizes. All this is beyond the scope of this paper.

\bibliographystyle{amsalpha}

\end{document}